\documentclass{amsart}
\title{Generalized Tambara-Yamagami categories}
\author{Jesse Liptrap}
\address{Mathematics Department\\
  University of California\\
  Santa Barbara, CA 93106\\
  U.S.A.}
\email{jliptrap@math.ucsb.edu}

\usepackage{enumerate,nicefrac}
\input{xy}
\xyoption{all}
\newtheorem{thmKludge}{Theorem}
\newtheorem{thm}{Theorem}[section]
\newtheorem*{thm*}{Theorem}
\newtheorem{prop}[thm]{Proposition}
\newtheorem*{prop*}{Proposition}
\newtheorem{cor}[thm]{Corollary}
\newtheorem*{cor*}{Corollary}
\newtheorem{lem}[thm]{Lemma}
\newtheorem*{lem*}{Lemma}
\theoremstyle{remark}

\newtheorem*{rem*}{Remark}

\newtheorem*{conv*}{Convention}
\newtheorem{exmp}[thm]{Example}
\newtheorem*{exmp*}{Example}
\theoremstyle{definition}
\newtheorem{defn}[thm]{Definition}
\newtheorem*{defn*}{Definition}
\newcommand{\abs}[1]{\lvert #1 \rvert}
\newcommand{\pf}{properly feudal}
\newcommand{\fe}{feudal}        % warning to macro writers: "improperly
\newcommand{\noast}{}
\newcommand{\comment}[1]{}
\newcommand{\N}{\ensuremath{\mathbb{N}}}
\newcommand{\Z}{\ensuremath{\mathbb{Z}}}
\newcommand{\F}{\ensuremath{\mathbb{F}}}
\newcommand{\ov}{\overline}

\newcommand{\f}[2]{f^{#1}_{#2}}
\DeclareMathOperator{\obj}{obj}
\DeclareMathOperator{\mor}{mor}
\DeclareMathOperator{\im}{im}
\DeclareMathOperator{\coker}{coker}
\newcommand{\ad}{\mathrm{ad}}
\newcommand{\lc}{\grave{\delta}}
\newcommand{\rc}{\acute{\delta}}
\newcommand{\ups}{\upsilon}
\newcommand{\om}{\omega}
\newcommand{\lexp}[2]{{\vphantom{#2}}^{#1}{#2}}
\newcommand{\lr}[3]{{}^{#1}{#2}^{#3}}
\newcommand{\setify}[1]{{\lfloor #1 \rfloor}}
\newcommand{\slashfrac}[2]{{{#1}/{#2}}}
\newcommand{\multip}[2]{{\langle #1, #2 \rangle}}
\newcommand{\maybreak}{\displaybreak[3]}
\newcommand{\hmm}[3]{{#1 \xrightarrow{#2} #3}}

\begin{document}

\begin{abstract}
  Fusion rules generalize groups by allowing multivalued multiplication.
  Groups are fusion rules of simple current index 1.  We classify
  nilpotent (in the sense of Gelaki and Nikshych) fusion rules of simple
  current index 2, and characterize the associated fusion categories.
\end{abstract}

\maketitle

\section*{Introduction}
\emph{Fusion categories} appear in conformal field theory, operator
algebras, representation theory, and quantum topology. They underly
models of fractional quantum Hall quasiparticles, the proposed raw
material for topological quantum computation.  \emph{Fusion rules} can
be viewed as nondeterministic groups.  Every fusion category has a
fusion rule up to isomorphism.  Thus the problem of classifying fusion
categories splits into two difficult subproblems: understand fusion
rules, and given a fusion rule, understand the associated fusion
categories---of which there are only finitely many up to equivalence
over an algebraically closed field of characteristic $0$ (Ocneanu
rigidity~\cite{ENO}).

To each fusion rule we assign a \emph{simple current index} so that
groups are precisely fusion rules of simple current index $1$.
Following Gelaki and Nikshych~\cite{GN} we assume \emph{nilpotence}.  We
classify---functorially---nilpotent fusion rules of simple current index
$2$, e.g., the \emph{Tambara-Yamagami} fusion rules~\cite{TY} and the
\emph{fermionic Moore-Read} fusion rule~\cite{Bo}, in terms of group
homomorphisms (Corollary~\ref{cor-ffrs}).  Generalizing Tambara and
Yamagamis' classification using nondegenerate symmetric bilinear forms
on groups~\cite{TY}, we give an elementary algebraic characterization of
the resulting fusion categories, and of $H^3(G,\F^\times)$ for $G$ a
group of even order and $\F$ a field (Corollaries~\ref{cor-grp},
\ref{cor-proper}).

For other generalizations of Tambara-Yamagami fusion categories, see
Siehler~\cite{S} and Etingof, Nikshych, and Ostrik \cite{ENOM}.

\section{Fusion rules}\label{sec:fusionRules}

\begin{defn}
  A \emph{multiset} over a set $L$ is a function $X \in \N^L$, assigning
  a \emph{multiplicity} $\multip{X}{x}$ to each $x \in L$.  A multiset
  is \emph{finite} if it has finite support.
\end{defn}

\begin{defn}
  A \emph{multimagma} is a set $L$ equipped with an operation $\ast
  \colon L \times L \to \N^L$.  If $X,Y \in \N^L$ are finite,
  define $X \ast Y \in \N^L$ as follows: for $z \in L$,
  \[\langle X \ast Y, z \rangle = \sum_{x,y \in L}
  \multip{X}{x} \multip{Y}{y} \multip{x \ast y}{z}.\] We say $L$ is
  \emph{locally finite} if $x \ast y$ is finite for $x,y \in L$.
\end{defn}

\begin{conv*}
  We identify each subset of a set $L$ with the multiset given by its
  indicator function, and each element $x \in L$ with the singleton
  $\{x\}$ given by the Kronecker delta $y \mapsto \delta_{x,y}$.  The
  symbol $\ast$ is often suppressed.
\end{conv*}

\begin{defn}\label{defn-fr}
  A \emph{fusion rule} is a locally finite multimagma satisfying
  \begin{align*}
    \forall x,y,z &\colon (x \noast y) \noast z = x \noast (y \noast z)\\
    \exists 1 \forall x  &\colon 1 \noast x  = x \noast 1 = x\\
    \forall x \exists \bar{x} \forall y &\colon \langle x \noast
    y,1\rangle = \langle y \noast x,1 \rangle = \delta_{y,\bar{x}}
  \end{align*}
\end{defn}

\begin{lem}\label{lem-basic} Let $L$ be a fusion rule. Then
  \begin{itemize}
  \item $xy \neq \emptyset$ for $x,y \in L$.
  \item $1 \in L$ is unique.
  \item If $x \in L$, then $\bar{x}$ is uniquely determined by $x$, and
    $\bar{\bar{x}}=x$.
  \end{itemize}
\end{lem}

\begin{defn}
  Define the \emph{underlying set} $\setify{X}$ of a multiset $X$ by $x
  \in \lfloor X \rfloor$ iff $\langle X,x\rangle > 0$.  Given a
  multimagma $L$, let $\lfloor L \rfloor$ be the multimagma with the
  same elements with operation $(x,y) \mapsto \lfloor xy \rfloor$.  We
  say $L$ is \emph{multiplicity-free} if $L = \lfloor L \rfloor$.
\end{defn}

\begin{rem*}
  Structures similar to multiplicity-free fusion rules, called
  \emph{hypergroups} or \emph{multigroups} \cite{DO, EO, CL}, have been
  studied since the 1930s, with various applications.
\end{rem*}

\begin{exmp}\label{exmp-ourFriends}  We will study three
  examples of multiplicity-free fusion rules:
  \begin{itemize}
  \item Groups are precisely fusion rules with single-valued fusion.
    Note $\bar{a} = a^{-1}$ for any group element $a$.
  \item Given a finite group $A$ and $m \notin A$, the
    \emph{Tambara-Yamagami} fusion rule has elements $A \cup \{m\}$
    fusing as follows: for $a,b \in A$,
    \begin{align*}
      a \ast b &= ab,& a \ast m &= m \ast a = m,& m \ast m &= A.
    \end{align*}
  \item The \emph{fermionic Moore-Read} fusion rule has six elements
    $\{\pm 1,\pm i,\pm i'\}$ fusing as follows: for $a,b \in \{\pm 1,
    \pm i\}$ and $p,q \in \{\pm i\}$,
    \begin{align*}
      a \ast b &= ab,& a \ast p' & = p' \ast a = (a^2p)',& p' \ast q' &=
      \{\pm \sqrt{pq}\}.
    \end{align*}
  \end{itemize}
\end{exmp}

\begin{defn}\label{defn-homom}
  Let $L$ and $L'$ be multiplicity-free multimagmas. A map $f \colon L
  \to L'$ is a \emph{homomorphism} if $f(xy) \subseteq f(x)f(y)$ for
  $x,y \in L$.
\end{defn}

\begin{defn}\label{defn-grading}
  Let $G$ be a group and $L$ a fusion rule.  A \emph{$G$-grading} on $L$
  is a surjective homomorphism $\lfloor L \rfloor \to G$.
\end{defn}

\begin{rem*}
  Many authors do not require gradings to be surjective.
\end{rem*}

\begin{defn} Let $L$ be a fusion rule and $S \subseteq L$.
  \begin{itemize}
  \item $S \ni 1$ is a \emph{subrule} if $z \in S$ whenever $\bar{z} \in
    S$ or $z \in \setify{xy}$ for some $x,y \in S$.
  \item A \emph{left coset} of $S$ in $L$ is a subset $\lfloor
    xS\rfloor$ for some $x \in L$.  The set $\nicefrac{L}{S}$ of all
    left cosets of $S$ in $L$ is a multimagma with operation $\ast_S$
    given by
    \[\langle X \ast_S Y , Z \rangle = \max \{\langle xy,z \rangle
    \mid x \in X,\ y \in Y,\ z \in Z \}.\]
  \item The \emph{index} of $S$ is the cardinality of $\nicefrac{L}{S}$.
  \end{itemize}
\end{defn}

\begin{defn}
  The \emph{adjoint subrule} $L_{\ad}$ of a fusion rule $L$ is the
  smallest subrule containing $\lfloor x\bar{x} \rfloor$ for all $x \in
  L$.  We say $L$ is \emph{nilpotent} if $L^{(n)} = \{1\}$ for some $n
  \in \N$, where $L^{(0)} = L$ and $L^{(m+1)} = (L^{(m)})_{\ad}$.
\end{defn}

\begin{thmKludge}[Dresher and Ore \cite{DO}, Gelaki and Nikshych
  \cite{GN}] \label{thmKludge-GN} Let $L$ be a fusion rule. Then
  $L_{\ad}$ is the intersection of all subrules $A$ such that $\lfloor
  \nicefrac{L}{A} \rfloor$ is a group.  The multimagma $\lfloor
  \nicefrac{L}{L_\ad} \rfloor$ is a group, called the \emph{universal
    grading group}, partitioning $L$.  Every grading of $L$ factors
  uniquely through the quotient projection $\setify{L} \to \lfloor
  \nicefrac{L}{L_\ad} \rfloor$, called the \emph{universal grading}.
\end{thmKludge}

\begin{defn}
  A fusion rule element $a$ is a \emph{simple current} if $a \bar{a} =
  1$.  The \emph{simple current index} of a fusion rule is the index of
  the set of simple currents.
\end{defn}

\begin{lem}\label{lem-maxsubgrp} Let $L$ be a fusion rule with simple
  currents $S$. Then
  \begin{enumerate}[{\normalfont (i)}]
  \item\label{it:simplecurrent} $a \in S$ iff $\bar{a}a=1$ iff $az$ and
    $za$ are singletons for all $z \in L$.
  \item\label{it:nomult} If $a \in S$, then $\langle xy,a\rangle \leq 1$
    for all $x,y \in L$.
  \item $S$ is the largest subrule of $L$ which is a group.
  \item $\nicefrac{L}{S}$ partitions $L$.
  \item $L$ is a group iff its simple current index is $1$.
  \end{enumerate}
\end{lem}

\section{Feudal fusion rules}

\begin{conv*}
  $\Z_2 = \{\pm 1\}$ is a group of order $2$.
\end{conv*}

\begin{defn}
  A fusion rule is \emph{\fe{}} if it is equipped with a $\Z_2$-grading
  $\gamma$ such that $\gamma^{-1}(1)$ is a group.  We call elements of
  $\gamma^{-1}(1)$ \emph{serfs} and elements of $\gamma^{-1}(-1)$
  \emph{lords}.  A fusion rule is \emph{\pf{}} if it is nilpotent with
  simple current index $2$.
\end{defn}

\begin{exmp}\label{exmp-things}\
  \begin{itemize}
  \item A $\Z_2$-graded group is improperly \fe{}.  Its adjoint subrule
    is trivial; it is its own simple current group and universal grading
    group.
  \item A Tambara-Yamagami fusion rule $A \cup \{m\}$ is \fe{} with
    serfs $A$ and lord $m$.  It is \pf{} iff $\abs{A} > 1$ iff $A$ is
    the simple current group.  The adjoint subrule is $A$; the universal
    grading group is $\{A, m\} \cong \Z_2$.  A \pf{} fusion rule with a
    lone lord is Tambara-Yamagami.
  \item The fermionic Moore-Read fusion rule is \fe{} with serfs $\{\pm
    1, \pm i\}$ and lords $\{\pm i'\}$.  It is \pf{}, with simple
    currents $\{\pm 1, \pm i\}$, adjoint subrule $\{\pm 1\}$, and
    universal grading group $\{\{\pm 1\}, i', \{\pm i\}, -i'\} \cong
    \Z_4$.
  \end{itemize}
\end{exmp}

\begin{lem}\label{lem-A}
  Let $L$ be a \fe{} fusion rule with serfs $S$. Then
  \begin{enumerate}[{\normalfont (i)}]
  \item\label{it:mult-free} $L$ is multiplicity-free.
  \item\label{it:trans} $S$ acts on lords by fusion, transitively on the
    left and on the right.
  \item\label{it:stab} $L_\ad \unlhd S$ is the left stabilizer and right
    stabilizer of any lord.
  \item\label{it:coset} Two lords $m,l$ fuse to a coset of $L_\ad$ in
    $S$, namely
    \begin{equation} \label{eq-stab} ml = \{a \in S \mid \bar{m}a = l\}
      = \{a \in S \mid m = a\bar{l}\}.
    \end{equation}
  \end{enumerate}
\end{lem}

\begin{proof}
  Let $M$ be the lords and $m,l \in M$.  By Lemma~\ref{lem-maxsubgrp},
  $S$ acts on $M$ on the left and on the right by fusion, and $ml$ is a
  subset of $S$, proving \eqref{it:mult-free}.  For $a \in S$,
  \begin{align*}
    a \in ml &\iff 1 \in \bar{a}ml \iff \bar{m}a = l \\
    &\iff 1 \in ml\bar{a} \iff m = a\bar{l},
  \end{align*}
  proving equation~\eqref{eq-stab}.  Then $ml \neq \emptyset$ implies
  \eqref{it:trans}.  Since $\{ml' \mid l' \in M\}$ and $\{m'l \mid m'
  \in M\}$ each partition $S$,
  \[ mb = m \iff mbl = ml \iff bl = l\] for $b \in S$, i.e., the right
  stabilizer of $m$ and the left stabilizer of $l$ coincide for
  arbitrary $m,l \in M$.  Let $A \unlhd S$ be this common stabilizer.
  Equation~\eqref{eq-stab} implies $m\bar{m} =A$ for all $m \in M$.
  Thus $A = L_\ad$, proving \eqref{it:stab}.  The orbit-stabilizer
  theorem of elementary group theory then completes \eqref{it:coset}.
\end{proof}

\begin{prop}
  A \pf{} fusion rule is uniquely \fe{}.  A \fe{} fusion rule is \pf{}
  or a $\Z_2$-graded group.
\end{prop}

\begin{proof}
  Let $L$ be \pf{} with simple currents $S$.  By
  Lemma~\ref{lem-maxsubgrp}, $M = {L \setminus S} \neq \emptyset$ and
  $am,ma \in M$ whenever $a \in S$ and $m \in M$.  Thus $\nicefrac{L}{S}
  = \{S,M\}$ with $S \ast_S S = S$ and $S \ast_S M = M \ast_S S = M$.

  To show $M \ast_S M = S$, we first need $L_\ad \subseteq S$.  Pick $m
  \in M$.  Since $M = mS$ and $ma\ov{ma} = m\bar{m}$ for $a \in S$, we
  see $L_\ad$ is the smallest subrule of $L$ containing $\lfloor
  m\bar{m}\rfloor$.  If $m \in L_\ad$, then $\lfloor m\bar{m} \rfloor
  \subseteq (L_\ad)_\ad$ implies $(L_\ad)_\ad = L_\ad$, contradicting
  nilpotence.  Therefore $L_\ad \subseteq S$.  Now pick any $m,l \in M$.
  Then $l = \bar{m}a$ for some $a \in S$, whence $\lfloor ml \rfloor =
  \lfloor m\bar{m}\rfloor a \subseteq S$.  Moreover $ml = \lfloor ml
  \rfloor$ by Lemma~\ref{lem-maxsubgrp}\eqref{it:nomult}.  Thus $M
  \ast_S M = S$, and $L \to \nicefrac{L}{S} \cong \Z_2$ is a \fe{}
  grading.

  Suppose $\gamma \colon L \to \Z_2$ is a different \fe{} grading.  Then
  $S' = \gamma^{-1}(1) \subset S$. Picking $a \in S \setminus S'$ and $m
  \in M$, we have $am \in M \cap S'$, a contradiction.  Thus $L$ is
  uniquely \fe{}.  Finally, a non-group \fe{} fusion rule is \pf{} by
  Lemma~\ref{lem-A}.
\end{proof}

\begin{defn} \label{defn-feudfun} Let $\mathcal{H}$ be the following
  category.  An object of $\mathcal{H}$ is a homomorphism
  $\hmm{S}{u}{G}$ of arbitrary groups $S,G$ such that $\abs{\coker{u}} =
  2$ and $\ker{u}$ is finite, with the innocuous technical conditions $S
  \cap (G \setminus \im{u}) = \emptyset$ and $\im{u} = S/ \ker{u}$.  A
  morphism from $\hmm{S}{u}{G}$ to
  $\hmm{\tilde{S}}{\tilde{u}}{\tilde{G}}$ in $\mathcal{H}$ is a pair of
  homomorphisms $(h_0, h_1)$ making the square
  \begin{equation}\label{diag-Hmorphism}
    \xymatrix{
      S \ar[r]^u \ar[d]_{h_0} & G \ar[d]^{h_1}\\
      \tilde{S} \ar[r]_{\tilde{u}} & \tilde{G}
    }
  \end{equation}
  commute, with $h_1(G \setminus \im{u}) \subseteq \tilde{G} \setminus
  \im{\tilde{u}}$.

  Let $\mathcal{L}$ be the category of \fe{} fusion rules and graded
  homomorphisms.  Let $\Phi \colon \mathcal{H} \to \mathcal{L}$ be the
  following functor.  For $H = (\hmm{S}{u}{G}) \in \obj{\mathcal{H}}$,
  let $\Phi H$ be the \fe{} fusion rule with serfs $S$ and lords $M = G
  \setminus \im{u}$ fusing as follows: for $a,b \in S$ and $m,l \in M$,
  \begin{align*}
    a \ast b &= ab,& a \ast m &= u(a)m,& m \ast a &= mu(a),& m \ast l &=
    u^{-1}(ml).
  \end{align*}
  For $\tilde{H} \in \obj \mathcal{H}$ and $(h_0,h_1) \in
  \mor_{\mathcal{H}}(H, \tilde{H})$, let $\Phi(h_0,h_1)$ agree with
  $h_0$ on $S$ and with $h_1$ on $M$.

  Inversely, let $\Gamma \colon \mathcal{L} \to \mathcal{H}$ be the
  following functor.  Given $L \in \obj \mathcal{L}$, let $\Gamma L =
  (\hmm{S}{u}{G})$ with $S$ the serfs (or the simple currents unless $L$
  is a group), $G$ the universal grading group, and $u$ the restriction
  to $S$ of the universal grading.  For $\tilde{L} \in \obj \mathcal{L}$
  and $t \in \mor_{\mathcal{L}}(L,\tilde{L})$, let $\Gamma t = (h_0,
  h_1)$ where $t$ agrees with $h_0$ and induces $h_1$.
\end{defn}

\begin{exmp}\label{ex-seqs}\
  \begin{enumerate}[{\normalfont (i)}]
  \item\label{it-grpseq} Let $G$ be a $\Z_2$-graded group with serfs
    $S$.  Then $\Gamma G$ is inclusion $S \to G$.
  \item Let $L = A \cup \{m\}$ be Tambara-Yamagami.  Then $\Gamma L$ is
    isomorphic to the trivial homomorphism $A \to \Z_2$.
  \item\label{it-MRseq} Let $L$ be the fermionic Moore-Read fusion rule.
    Then $\Gamma L$ is isomorphic to the nontrivial nonidentity
    homomorphism $\Z_4 \to \Z_4$.
  \end{enumerate}
\end{exmp}

\begin{thmKludge}\label{thmKludge-feudalStructure}
  The category of feudal fusion rules and graded homomorphisms is
  isomorphic to the category $\mathcal{H}$ of
  Definition~\ref{defn-feudfun}, via the functors therein.
\end{thmKludge}

\begin{cor}\label{cor-ffrs}
  Up to isomorphism, properly \fe{} fusion rules are in 1-1
  correspondence with group homomorphisms whose cokernels have order $2$
  and whose kernels are nontrivial and finite.
\end{cor}

\begin{proof} [Proof of Theorem~\ref{thmKludge-feudalStructure}]
  First we check $\Phi \colon \mathcal{H} \to \mathcal{L}$ is a functor.
  Let $H = (\hmm{S}{u}{G}) \in \obj \mathcal{H}$.  We check $\Phi H \in
  \obj \mathcal{L}$.  Let $A = \ker{u}$ and $M = G \setminus \im{u}$.
  For $a,b,c \in S$ and $m,l,r \in M$,
  \begin{align*}
    (a \ast b) \ast c = abc = a \ast (b \ast c),\\
    (a \ast m) \ast b = u(a)mu(b) = a \ast (m \ast b),\\
    (a \ast b) \ast m =  u(ab)m = a \ast (b \ast m),\\
    (m \ast a) \ast b = m \ast (a \ast b),\maybreak{}\\
    (m \ast a) \ast l = u^{-1}(m u(a) l) = m \ast (a \ast l),\\
    (a \ast m) \ast l = u^{-1}(u(a)ml) = au^{-1}(ml) = a \ast (m \ast
    l),\\
    (m \ast l) \ast a = m \ast (l \ast
    a),\\
    (m \ast l) \ast r = (g \mapsto \abs{A}\delta_{g, mlr}) = m \ast (l
    \ast r).
  \end{align*}
  Thus $\ast$ is associative.  Therefore $\Phi H \in \obj \mathcal{L}$.
  Now let $\tilde{H} = (\hmm{\tilde{S}}{\tilde{u}}{\tilde{G}}) \in
  \obj{\mathcal{H}}$ and $h = (h_0, h_1) \in \mor_\mathcal{H}(H,
  \tilde{H})$.  For $a,b \in S$ and $m,l \in M$,
  \begin{align*}
    \Phi h(a \ast b) = h_0(ab) = \Phi h(a) \ast \Phi h(b),\\
    \Phi h(a \ast m) = h_1(u(a)m) = \tilde{u}(h_0(a))h_1(m) = \Phi h(a) \ast \Phi h(m),\\
    \Phi h(m \ast a) = \Phi h(m) \ast \Phi h(a).
  \end{align*}
  Since $\abs{\coker{u}} = 2$, there exists $b \in u^{-1}(ml)$.  Let
  $\tilde{A} = \ker \tilde{u}$.  Since square~\eqref{diag-Hmorphism}
  commutes, $h_0(A) \subseteq \tilde{A}$.  Then
  \begin{align*}
    \Phi h(m \ast l) = h_0(b A) \subseteq h_0(b) \tilde{A} =
    \tilde{u}^{-1}(\tilde{u}(h_0(b))) = \tilde{u}^{-1}(h_1(ml)) = \Phi
    h(m) \ast \Phi h(l)
  \end{align*}
  Thus $\Phi h$ is a homomorphism.  Therefore $\Phi \colon \mathcal{H}
  \to \mathcal{L}$ is a functor.

  Now we check $\Gamma \colon \mathcal{L} \to \mathcal{H}$ is a functor.
  Suppose $L,\tilde{L} \in \obj \mathcal{L}$, with serfs $S, \tilde{S}$,
  lords $M, \tilde{M}$, adjoint subrules $A, \tilde{A}$, and
  restrictions $u, \tilde{u}$ to serfs of the universal gradings,
  respectively.  Then $\Gamma L = (\hmm{S}{u}{G})$ and $\Gamma \tilde{L}
  = (\hmm{\tilde{S}}{\tilde{u}}{\tilde{G}})$ are in $\obj{\mathcal{H}}$.
  Suppose $t \in \mor_{\mathcal{L}}(L,\tilde{L})$.  Since $t$ is graded,
  it restricts to a homomorphism $h_0 \colon S \to \tilde{S}$.  By
  Lemma~\ref{lem-A}\eqref{it:stab}, $h_0(A) \subseteq \tilde{A}$.  Let
  $h_{0.5}$ be the induced homomorphism $\nicefrac{S}{A} \to
  \nicefrac{\tilde{S}}{\tilde{A}}$.  Recalling $G = (\nicefrac{S}{A})
  \cup M$ and $\tilde{G} = (\nicefrac{\tilde{S}}{\tilde{A}}) \cup
  \tilde{M}$, let $h_1 \colon G \to \tilde{G}$ agree with $h_{0.5}$ on
  $\nicefrac{S}{A}$ and with $t$ on $M$.  Then $\Gamma t = (h_0, h_1)$
  is defined and square~\eqref{diag-Hmorphism} commutes.  To check $h_1$
  is a homomorphism, let $a,b \in S$ and $m,l \in M$ and $c \in m \ast
  l$.  By Lemma~\ref{lem-A},
  \begin{align*}
    h_1((aA)(bA)) = h_{0.5}((aA)(bA)) = h_1(aA)h_1(bA),\\
    h_1((aA)m) = t(a \ast m) = t(a) \ast t(m) = h_1(aA)h_1(m),\\
    h_1(m(aA)) = h_1(m)h_1(aA),\\
    h_1(ml) = h_{0.5}(cA) = t(c)\tilde{A} = t(m) \ast t(l) =
    h_1(m)h_1(l).
  \end{align*}
  Therefore $\Gamma$ is a functor.  It is easy to see $\Phi$ and
  $\Gamma$ are mutually inverse.
\end{proof}

\section{Fusion systems}

\begin{conv*}
  $\F$ is a field.
\end{conv*}

\begin{defn}
  A \emph{fusion category} is a rigid semisimple $\F$-linear monoidal
  category with simple monoidal unit and one-dimensional endomorphism
  spaces of simple objects.  Equivalence of fusion categories is
  $\F$-linear monoidal equivalence.
\end{defn}

\begin{rem*}
  Fusion categories are normally taken finite, i.e., having finitely
  many isomorphism classes of simple objects.  Ocneanu rigidity assumes
  finiteness \cite{ENO}.  In this paper finiteness is not required.
\end{rem*}

\begin{defn}
  Let $L$ be a complete set of representatives of isomorphism classes of
  simple objects of a fusion category $C$.  For $x,y \in L$ define $xy
  \in \N^L$ by $\langle xy,z \rangle = \dim_{\F} \mathrm{mor}_C(z,x
  \otimes y)$ where $\otimes$ is the monoidal product on $C$.  Then $L$
  is a fusion rule, unique up to isomorphism.  We call $C$ a fusion
  category \emph{on} $L$.
\end{defn}

  \begin{rem*}
    Although fusion categories are the main topic of interest, we
    circumvent them via a technical device of Yamagami.  For notational
    simplicity we only treat the multiplicity-free case, since \fe{}
    fusion rules are multiplicity-free.
  \end{rem*}

\begin{defn}\label{defn-fs}  
  Let $L$ be a multiplicity-free fusion rule.  A \emph{fusion system} on
  $L$ is map $f \colon L^6 \to \F$, assigning a \emph{recoupling
    coefficient} $\f{xyz}{urv}$ to each sextuple $(x,y,z,u,r,v)$, such
  that for $w,x,y,z,p,u,r,v,q \in L$,
  \begin{description}
  \item[Admissibility] $\f{xyz}{urv} = 0$ unless $(x,y,z,u,r,v)$ is
    \emph{admissible}, i.e., $u \in xy$ and $v \in yz$ and $r \in uz
    \cap xv$.
  \item[Invertibility] Each \emph{recoupling matrix} $F^{xyz}_r =
    (\f{xyz}{urv})_{v,u}$ is invertible, where $v,u$ range over all
    elements making $(x,y,z,u,r,v)$ admissible.
  \item[Pentagon axiom] $P^{wxyz}_{purvq}:\qquad \f{wxq}{prv}
    \f{pyz}{urq} = \sum_{s \in L} \f{xyz}{svq} \f{wsz}{urv}
    \f{wxy}{pus}$
  \item[Triangle axiom] $F^{x1y}_r$ is the identity matrix, $1 \times 1$
    or $0 \times 0$.
  \item[Rigidity] $f^{r\bar{r}r}_{1r1} =
    ((F^{\bar{r}r\bar{r}}_{\bar{r}})^{-1})_{11} \neq 0$.
  \end{description}
\end{defn}

\begin{defn}\label{defn-equiv}
  Let $f$ and $\tilde{f}$ be fusion systems on a multiplicity-free
  fusion rule $L$. A \emph{gauge transformation} from $f$ to
  $\tilde{f}$, rendering them \emph{gauge equivalent}, is a map $\xi
  \colon L^3 \to \F$, written $(x,y,r) \mapsto \xi^{xy}_r$, such that
  for $x,y,z,u,r,v \in L$,
  \begin{description}
  \item[Invertibility] $\xi^{xy}_r \neq 0$ iff $r \in xy$.
  \item[Rectangle axiom] $G^{xyz}_{urv}: \qquad \f{xyz}{urv} \xi^{yz}_v
    \xi^{xv}_r = \xi^{xy}_u \xi^{uz}_r \tilde{f}^{xyz}_{urv}$
  \item[Normalization] $\xi^{1r}_r = \xi^{r1}_r = 1$.
  \end{description}
  We say $f$ and $\tilde{f}$ are \emph{equivalent} if they are gauge
  equivalent up to relabeling by an automorphism of $L$.
\end{defn}

\begin{rem*}
  Yamagami's polygonal notation \cite{Y} is a good way to visualize the
  pentagon and rectangle axioms.
\end{rem*}

\begin{thmKludge}[Yamagami \cite{Y}]\label{thmKludge-Y}
  Modulo equivalence, fusion categories on a multi\-plicity-free fusion
  rule $L$ are in 1-1 correspondence with fusion systems on $L$.
\end{thmKludge}

\begin{proof}
  This is Proposition~1.1 and Lemma~2.2 of \cite{Y} written out in
  coordinates with the added assumption that $L$ is multiplicity-free.
\end{proof}

\begin{lem}\label{lem-1top}
  $F^{1xy}_r = F^{xy1}_r = \mathrm{Id}$ for $F$ as in
  Definition~\ref{defn-fs} and $x,y,r \in L$.
\end{lem}

\begin{proof}
  $P^{\bar{r},1,x,y}_{\bar{r},\bar{y},1,r,r}$ and
  $P^{x,y,1,\bar{r}}_{r,r,1,\bar{x},\bar{r}}$ and the triangle axiom.
\end{proof}

\begin{defn}\label{defn-cohom}\label{defn-lrcob}
  Let $S$ be a group, $U$ an $S$-bimodule, and $n \in \{1,2,3\}$.  The
  \emph{left and right coboundary operators} $\lc,\rc \colon U^{S^n}
  \to U^{S^{n+1}}$ are defined as follows: for an \emph{$n$-cochain}
  $h_n \colon S^n \to U$, define $\lc h_n, \rc h_n \colon S^{n+1} \to
  U$ by
  \begin{align*}
    \lc h_1(a,b) &= \frac{ h_1(a) {}^ah_1(b)}{h_1(ab)} &  \rc h_1(a,b) &= \frac{ h_1^b(a) h_1(b)}{h_1(ab)}\\
    \lc h_2(a,b,c) &= \frac{ h_2(a,bc){}^ah_2(b,c)}{h_2(a,b)h_2(ab,c)} &\rc h_2(a,b,c) &= \frac{ h_2(a,bc)h_2(b,c)}{h_2^c(a,b)h_2(ab,c)}\\
    \lc h_3(a,b,c,d) &= \frac{ h_3(a,b,c) h_3(a,bc,d) {}^ah_3(b,c,d)} {
      h_3(a,b,cd)h_3(ab,c,d)} & \rc h_3(a,b,c,d) &= \cdots
  \end{align*}
  Here $U$ is written multiplicatively; left and right exponentiation of
  cochains denotes the $S$-actions on $U$.  If it is known that $\lc h_n
  = \rc h_n$, we write $\delta h_n = \lc h_n$.

  $H^n(S,U)$ is the abelian group $\slashfrac{\ker (\lc^n)}{\mathrm{im}
    (\lc^{n-1})}$.

  A \emph{normalized} cochain is $1$ whenever any argument is $1$.
\end{defn}

\begin{exmp}[well-known]\label{ex-grpcohom}
  View $\F^\times$ as a trivial module over a group $G$.  Then fusion
  systems on $G$ up to gauge equivalence are in 1-1 correspondence with
  $H^3(G,\F^\times)$, and fusion categories (or systems) on $G$ up to
  equivalence are in 1-1 correspondence with
  $H^3(G,\F^\times)/{\mathrm{aut}G}$.
\end{exmp}

\begin{proof}
  Recall an $n$-cochain $h$ is an $n$-\emph{cocycle} if $\lc h = 1$ and
  an $n$-\emph{coboundary} if $h = \lc k$ for some $k$, and cocycles $h,
  \tilde{h}$ are \emph{cohomologous} if $\nicefrac{h}{\tilde{h}}$ is a coboundary.
  Identify maps $f \colon G^6 \to \F^\times$ satisfying the
  admissibility axiom of Definition~\ref{defn-fs} with $3$-cochains, and
  maps $\xi \colon G^3 \to \F^\times$ satisfying the invertibility axiom
  of Definition~\ref{defn-equiv} with $2$-cochains, via
  \begin{align*}
    f^{a,b,c}_{u,r,v} &= \delta_{u,ab}\delta_{r,abc}\delta_{v,bc}f(a,b,c), &
    \xi^{a,b}_r &= \delta_{r,ab} \xi(a,b).
  \end{align*}
  If $f$ is a normalized $3$-cocycle, $\delta f(r,\bar{r},r,\bar{r}) =
  1$ implies $f(r,\bar{r},r)f(\bar{r},r,\bar{r}) = 1$, i.e., $f$ is
  rigid.  Therefore every normalized $3$-cocycle is a fusion system on
  $G$, and conversely by Lemma~\ref{lem-1top}.  Moreover $\xi$ is a
  gauge transformation from $f$ to $\tilde{f}$ iff $\xi$ is normalized
  and $\tilde{f} = f \delta \xi$.  By Lemma 15.7.1 of~\cite{H}, every
  cocycle is cohomologous to a normalized cocycle; by Lemma 15.7.2
  of~\cite{H}, every normalized coboundary is the coboundary of a
  normalized cochain.  Our identification of fusion systems on $G$ with
  normalized $3$-cocycles thus descends to a 1-1 correspondence between
  fusion systems on $G$ up to gauge equivalence and $H^3(G,\F^\times)$,
  which descends to the second claimed correspondence via
  Theorem~\ref{thmKludge-Y}.
\end{proof}

\begin{lem}\label{lem-gaugify}
  If $f$ is a fusion system and $\xi$ satisfies the invertibility and
  normalization axioms of a gauge transformation, then $\tilde{f}$
  defined by the rectangle axiom is a fusion system, and $\xi$ is a
  gauge transformation from $f$ to $\tilde{f}$.
\end{lem}

\begin{proof}
  Routine.
\end{proof}

\section{Feudal fusion systems}

\begin{defn}
  A \emph{symmetric bicharacter} on a finite group $A$ over a ring $B$
  is a map $\chi \colon A \times A \to B^\times$ such that
  \begin{align*}
    \chi(b,a) &= \chi(a,b), & \chi(ab,c) &= \chi(a,c)\chi(b,c).
  \end{align*}
  We say $\chi$ is \emph{nondegenerate} if $\sum_b \chi(a,b) = 0$ for $a
  \neq 1$.
\end{defn}

\begin{defn}\label{defn-bimod}
  An \emph{involutory ambidextrous algebra} over a group $S$ is a ring
  $B$ with two operations
  \begin{align*}
    B \to B &\colon \mu \mapsto \bar{\mu}\\
    S \times B \times S \to B &\colon (a,\mu,b) \mapsto {}^a\mu^b
  \end{align*}
  such that $\mu \mapsto \bar{\mu}$ is an involution (ring
  antiautomorphism of order two), $\mu \mapsto {}^a\mu^b$ is a ring
  endomorphism for $a,b \in S$, and
  \begin{align*} {}^a({}^b\mu{}^c)^d &= {}^{ab}\mu^{cd},& \ov{{}^a\mu^b}
    &= {}^{\bar{b}}\bar{\mu}^{\bar{a}}
  \end{align*}
  for $a,b,c,d \in S$ and $\mu \in B$.
\end{defn}

\begin{conv*}
  Let $B$ be an involutory ambidextrous algebra over a group $S$, and
  let $X$ be a set.  Then $B^X$ inherits the involutory ambidextrous
  $S$-algebra structure of $B$: for $\epsilon \in B^X$ and $a,b \in S$,
  define ${}^a\epsilon^b, \bar{\epsilon} \in B^X$ by ${}^a\epsilon^b(x)
  = {}^a\epsilon(x)^b$ and $\bar{\epsilon}(x) = \ov{\epsilon(x)}$ for $x
  \in X$.  For $\mu,\nu \in B$ and $a \in S$, we write $\mu {}^a \nu$
  for $\mu({}^a\nu)$, not $(\mu^a)\nu$.
\end{conv*}

\begin{defn}
  Let $B$ be an involutory ambidextrous algebra over a group $S$.  Let
  $\chi,\ups \colon S \times S \to B^\times$ and $\tau \in B^\times$.
  \begin{itemize}
  \item $\chi$ is a \emph{$\tau$-quasisymmetric $\ups$-biderivation} if
    $\ups$ is normalized and for $a,b,c \in S$,
    \begin{align*} \bar{\chi}(b,a) &= {}^{\bar{a}}\chi^{\bar{b}}(a,b)
      \frac{{}^{\bar{a}}\tau^{\bar{b}} \cdot \tau}{{}^{\bar{a}}\tau
        \tau^{\bar{b}}}&
      \frac{\ups}{\ups^c}(a,b) \chi(ab,c) &= \chi(a,c) {}^a\chi(b,c).
    \end{align*}
  \item Suppose the set $A$ of elements of $S$ acting trivially on $B$
    is finite.  The triple $(\chi,\ups,\tau)$ is an
    \emph{\"uberderivation} on $S$ over $B$ if $\chi$ is a
    $\tau$-quasisymmetric $\ups$-bideri\-vation such that the symmetric
    bicharacter $\chi|_{A \times A}$ is nondegenerate, and $\abs{A} \tau
    \bar{\tau} = 1_B$.
  \item Let $\mathrm{fix}(S)$ be the elements of $B$ fixed under the
    $S$-actions.  A \emph{gauge transformation} from $(\chi,\ups,\tau)$
    to another \"uberderivation
    $(\tilde{\chi},\tilde{\ups},\tilde{\tau})$, rendering them
    \emph{gauge equivalent}, is a triple $(\theta,\phi,\varsigma) \in
    {\mathrm{fix}(S)}^{S \times S} \times (B^\times)^S \times
    (B^\times)$, with $\phi$ normalized, such that for $a,b \in S$,
    \begin{align*}
      \frac{\tilde{\chi}(a,b)}{\chi(a,b)} &= \frac{ \phi(a)
        \lr{a}{\bar{\phi}}{b}(b) \lexp{a}{\varsigma}^b \cdot \varsigma
      }{ \phi^b(a) \bar{\phi}^b(b) \lexp{a}{\varsigma} \varsigma^b} &
      \frac{\tilde{\ups} }{\ups} &= \frac{\lc \phi}{\theta} &
      \frac{\tilde{\tau}}{\tau} &= \frac{\bar{\varsigma}}{\varsigma}
    \end{align*}
  \end{itemize}
\end{defn}

\begin{defn}
  Let $L$ be a feudal fusion rule with serfs $S$ and lords $M$.  Let $B
  = \F^M$, with ring structure inherited from $\F$.  Then $B$ is an
  involutory ambidextrous $S$-algebra: for $\mu \in B$ and $a,b \in S$,
  define ${}^a\mu^b, \bar{\mu} \in B$ by ${}^a\mu^b(m) =
  \mu(\bar{a}m\bar{b})$ and $\bar{\mu}(m) = \mu(\bar{m})$ for $m \in M$.
  For a fusion system $f$ on $L$ or a gauge transformation $\xi$ of
  fusion systems on $L$, let $\Psi f = (\chi,\ups,\tau)$ or $\Psi\xi =
  (\theta,\phi,\varsigma)$, where for $a,b \in S$ and $m \in M$,
  \begin{align*}
    \chi(a,b)(m) &= f^{a,\bar{a}m\bar{b},b}_{m\bar{b},m,\bar{a}m} &
    \theta(a,b)(m) &= \xi^{a,b}_{ab}\\
    \ups(a,b)(m) &= f^{a,b,\bar{b}\bar{a}m}_{ab,m,\bar{a}m} &
    \phi(a)(m) &= \xi^{a,\bar{a}m}_m\\
    \tau(m) &= f^{m,\bar{m},m}_{1,m,1} & \varsigma(m)
    &=\xi^{m,\bar{m}}_1
  \end{align*}
\end{defn}

\begin{thmKludge}\label{thmKludge-main}
  Let $\mathcal{Y}$ be the category of fusion systems on a feudal fusion rule
  with lords $M$, and $\mathcal{X}$ the category of \"uberderivations on
  serfs over $\F^M\!$; in each category morphisms are gauge
  transformations, composed via multiplication in $\F$.  Then $\Psi
  \colon \mathcal{Y} \to \mathcal{X}$ is an equivalence, surjective on the nose.
\end{thmKludge}

\begin{cor}\label{cor-proper}
  Let $L$ be a \pf{} fusion rule, or a $\Z_2$-graded group all of whose
  automorphisms are graded, with serfs $S$ and lords $M$.  Then fusion
  categories \textup{(}or fusion systems\textup{)} on $L$ up to
  equivalence are in 1-1 correspondence with \"uberderivations on $S$
  over $\F^M$ up to gauge equivalence and up to simultaneously permuting
  $S$ and $M$ by an automorphism of $L$.
\end{cor}

\begin{proof}
  Theorems \ref{thmKludge-main} and \ref{thmKludge-Y}.
\end{proof}

\begin{exmp}[Tambara-Yamagami \cite{TY}]\label{exmp-TY}
  Let $L = A \cup \{m\}$ be a Tambara-Yamagami fusion rule.
  Then fusion categories on $L$ up to equivalence are in 1-1
  correspondence with pairs $(\chi,\tau)$ up to relabeling $\chi$ by an
  automorphism of $A$, where $\chi \colon A \times A \to \F^\times$ is a
  nondegenerate symmetric bicharacter and $\tau = \pm \abs{A}^{-1/2} \in
  \F$.
\end{exmp}

\begin{proof}
  This is a special case of Corollary~\ref{cor-proper}.
\end{proof}

\begin{exmp}[Bonderson~\cite{Bo}] \label{exmp-MR} If $\F$ has square
  roots, fermionic Moore-Read fusion categories up to equivalence are in
  1-1 correspondence with $4$th roots of $-1$ in $\F$.
\end{exmp}

\begin{cor}\label{cor-grp}
  Let $G$ be a group and $S$ an index $2$ subgroup.  Then
  $H^3(G,\F^\times)$ is in 1-1 correspondence with \"uberderivations on
  $S$ over $\F^{G\setminus S}$ up to gauge equivalence.
\end{cor}

\begin{proof}
  Theorem~\ref{thmKludge-main} and Example~\ref{ex-grpcohom}.
\end{proof}

\begin{exmp}[well-known] \label{exmp-Z4} If $\F$ has square roots,
  $H^3(\Z_4,\F^\times)$ is in 1-1 correspondence with $4$th roots of
  unity in $\F$.
\end{exmp}

\begin{cor}\label{cor-ab}
  If there is a fusion category on a feudal fusion rule with adjoint
  subrule $A$ and lords $M$, then
  \begin{itemize}
  \item $A$ is abelian.
  \item The characteristic of $\F$ does not divide $\abs{A}$.
  \item If $|M|$ is odd or $m = \bar{m}$ for some $m \in M$, then
    $\sqrt{\abs{A}} \in \F$.
  \end{itemize}
\end{cor}

\begin{proof}
  Let $(\chi,\ups,\tau)$ be an \"uberderivation on serfs over $B =
  \F^M$.  The two desired conditions on $\abs{A}$ follow from $\abs{A}
  \tau \bar{\tau} = 1_B$.  Since $\chi|_{A \times A}$ is a bicharacter,
  \[ \sum_{c \in A}\chi(ab\bar{a}\bar{b},c) = \abs{A}1_{B}\] for $a,b
  \in A$.  Since $\abs{A} \neq 0$ in $\F$ and $\chi|_{A \times A}$ is
  nondegenerate, $ab\bar{a}\bar{b}=1$.
\end{proof}

\begin{rem*}
  The results of this section would still hold if we did not require
  fusion categories and fusion systems to be rigid.
\end{rem*}

\begin{lem}\label{lem-consttau}
  Consider a feudal fusion rule with two lords $M = \{m_1,m_2\}$.
  Suppose $\F$ has square roots, and let $\tau_0(m_1) = \tau_0(m_2) =
  \pm \abs{A}^{-1/2} \in \F^\times$.  Then any \"uberderivation on serfs
  over $\F^M$ is gauge equivalent to one of the form $(\chi, \upsilon,
  \tau_0)$.  If two such \"uberderivations are gauge equivalent, they
  are related by a gauge transformation of the form $(\theta, \phi, 1)$.
\end{lem}

\begin{proof}
  Let $(\chi, \upsilon, \tau) \in \obj \mathcal{X}$, where $\mathcal{X}$ is
  as in Theorem~\ref{thmKludge-main}.  Choose $\varsigma \in \F^M$ such
  that $\nicefrac{\varsigma(m_1)}{\varsigma(m_2)} = \pm
  \sqrt{\nicefrac{\tau(m_1)}{\tau(m_2)}}$.  Let $\tilde{\tau}(m) =
  \nicefrac{\varsigma(\bar{m}) \tau(m)}{\varsigma(m)}$ for $m \in M$.
  Then $\tilde{\tau}$ is constant on $M$.  By Lemma~\ref{lem-gaugify}
  and Theorem~\ref{thmKludge-main} there exist $\theta, \phi, \varsigma,
  \tilde{\chi}, \tilde{\ups}$ such that $(\tilde{\chi},
  \tilde{\upsilon}, \tilde{\tau}) \in \obj \mathcal{X}$ and $(\theta,
  \phi, \varsigma) \in \mor_{\mathcal{X}}((\chi, \upsilon, \tau),
  (\tilde{\chi}, \tilde{\upsilon}, \tilde{\tau}))$.  Since
  $\tilde{\tau}^2 \equiv \nicefrac{1}{\abs{A}}$ and
  $\nicefrac{\bar{\varsigma}}{\varsigma}$ has sign freedom,
  w.l.o.g. $\tilde{\tau} = \tau_0$.

  Now let $\chi, \upsilon, \tilde{\chi}, \tilde{\upsilon}$ be arbitrary
  such that $(\chi, \upsilon, \tau_0), (\tilde{\chi}, \tilde{\upsilon},
  \tau_0) \in \obj{\mathcal{X}}$ are related by some gauge equivalence
  $(\theta, \phi, \varsigma)$.  Then $\bar{\varsigma} = \varsigma$,
  implying
  \[
  \frac{\tilde{\chi}(a,b)}{\chi(a,b)} = \frac{ \phi(a)
    {}^a\bar{\phi}^b(b) }{ \phi^b(a) \bar{\phi}^b(b)}
  \]
  Thus $(\theta, \phi, 1) \in \mor_{\mathcal{X}}((\chi, \upsilon, \tau),
  (\tilde{\chi}, \tilde{\upsilon}, \tilde{\tau}))$.
\end{proof}

\begin{conv*}
  When $\abs{M} = 2$, identify $\F^M$ with $\F^2$ naturally, and $\F$
  with the diagonal.
\end{conv*}

\begin{proof}[Proof of Example~\ref{exmp-MR}]
  Let $S = \{1, -1, i, -i\}$ and $M = \{\pm i'\}$.  Suppose
  $(\chi,\ups,\tau)$ is an \"uberderivation on $S$ over $B = \F^M$.  By
  Lemma~\ref{lem-consttau} we may take $\tau$ constant.  Writing $\chi$
  and $\nicefrac{\ups}{\bar{\ups}}$ as matrices over $B$ indexed by $S$,
  we find
  \begin{align*}
    \chi &= \begin{pmatrix}
      1 & 1 & 1 & 1\\
      1 & -1 & p & -p\\
      1 & p & x & r\\
      1 & -p & \bar{r} & x
    \end{pmatrix} & \nicefrac{\ups}{\bar{\ups}} &= \begin{pmatrix}
      1 & 1 & 1 & 1\\
      1 & p^2 & -\nicefrac{pr}{x} & -\nicefrac{px}{r}\\
      1 & -\nicefrac{\bar{p}r}{x} & \nicefrac{x^2}{p} & xr\\
      1 & -\nicefrac{\bar{p}x}{r} & x\bar{r} & -\nicefrac{x^2}{p}
    \end{pmatrix}
  \end{align*}
  for some $x \in \F^\times$ and $p,r \in B^\times$; the only
  requirements are $x^4 = p\bar{p} = -1$ and $r\bar{r} = -x^2$.  By
  Lemma~\ref{lem-gaugify} and Theorem~\ref{thmKludge-main} the gauge
  equivalence class of $f$ is uniquely determined by $x$, which is
  invariant under fusion rule automorphisms.  Then invoke
  Corollary~\ref{cor-proper}.
\end{proof}

\begin{proof}[Proof of Example~\ref{exmp-Z4}]
  Regard $\Z_4 = \{\pm 1, \pm i\}$ as \fe{} with lords $M = \{\pm i\}$.
  Then an \"uberderivation $(\chi,\ups, \tau)$ on $\{\pm 1\}$ over
  $\F^M$ is uniquely determined by $p,q, \tau \in \F^M$, where $p =
  \chi(-1,-1)$ and $q = \ups(-1,-1)$.  By Lemma~\ref{lem-consttau} we
  may take $\tau$ constant.  Then the only requirements on $p,q$ are $p
  \in \F^\times$ and $p^2 = \nicefrac{q}{\bar{q}}$, i.e., $p^4=1$ and
  $q_1 = p^2q_2$ where $(q_1,q_2) = q$.  The gauge equivalence class of
  $(\chi, \upsilon, \tau)$ is uniquely determined by $p$.  Then invoke
  Corollary~\ref{cor-grp}.
\end{proof}

\section{Proof of Theorem~\ref{thmKludge-main}}
We generalize Tambara and Yamagamis' proof \cite{TY} of
Example~\ref{exmp-TY}.

\begin{conv*}
  We identify $\F^\times$ with the set of constant functions $M \to
  \F^\times$.  Lemma~\ref{lem-A} is used throughout.
\end{conv*}

Let $L$ be a feudal fusion rule with serfs $S$, lords $M$, and adjoint
subrule $A$, and let $B = \F^M$.  We write $f,\tilde{f}$ for arbitrary
fusion systems on $L$ with recoupling matrices $F, \tilde{F}$;
$a,b,c,d,e$ for arbitrary serfs; and $m,l$ for arbitrary lords.  We
identify $f$ (likewise $\tilde{f}$) with the collection of eight
functions
\begin{align*}
  \alpha &\colon S \times S \times S \to \F^\times \\
  \alpha_1,\alpha_2,\alpha_3,\beta_1,\beta_2,\beta_3 &\colon S \times S \to B^\times\\
  \gamma &\colon S \times S \to B
\end{align*}
defined according to the following convention:
\begin{align*}
  & &\alpha(a,b,c) &= F^{a,b,c}_{abc}&&\\
  \alpha_1(a,b)(m) &= F^{m\bar{b}\bar{a},a,b}_m  & \alpha_2(a,b)(m) &= F^{a,\bar{a}m\bar{b},b}_m & \alpha_3(a,b)(m) &= F^{a,b,\bar{b}\bar{a}m}_m\\
  \beta_1(a,b)(m) &= F^{a,m,\bar{m}\bar{a}b}_b & \beta_2(a,b)(m) &= F^{m,a,\bar{a}\bar{m}b}_b & \beta_3(a,b)(m) &= F^{b\bar{a}\bar{m},m,a}_b\\
  &&\gamma(a,b)(m) &= f^{m\bar{a},a\bar{m}b,\bar{b}m}_{b,m,a}&&
\end{align*}
By Lemma~\ref{lem-1top}, $\alpha, \alpha_1,\alpha_2,\alpha_3$ are
normalized and $\beta_1(1,-) = \beta_2(1,-) = \beta_3(1,-) \equiv 1$.

\begin{lem}\label{lem:cosetmats}
  Every nonempty recoupling matrix is $1 \times 1$ or $(\gamma(a,b))_{a
    \in A', b \in A''}$ for some cosets $A',A''$ of $A$ in $S$.
\end{lem}

We write $\xi$ for an arbitrary gauge transformation from $f$ to
$\tilde{f}$, identified with the collection of four functions
\begin{align*}
  \theta &\colon S \times S \to \F^\times, & \phi,\psi,\omega&\colon S
  \to B^\times
\end{align*}
defined according to the following convention:
\begin{align*}
  \theta(a,b)(m) &= \xi^{a,b}_{ab} & \omega(a)(m) &=
  \xi^{m,\bar{m}a}_a \\
  \phi(a)(m) &= \xi^{a,\bar{a}m}_m & \psi(a)(m) &= \xi^{m\bar{a},a}_m
\end{align*}
By Definition~\ref{defn-equiv}, $\theta,\phi,\psi$ are normalized.  It
is routine to gather the rectangle axiom for $\xi$ and the pentagon
axiom for $f$ into
\begin{align*} G^{000} &= G^{a,b,c}_{ab,abc,bc}: & \tilde{\alpha} &=
  \alpha \delta \theta\\ G^{100} &=
  G^{m\bar{b}\bar{a},a,b}_{m\bar{b},m,ab}: & \tilde{\alpha}_1 \rc \psi
  &= \theta \alpha_1\\ G^{010} &=
  G^{a,\bar{a}m\bar{b},b}_{m\bar{b},m,\bar{a}m}: &
  \tilde{\alpha}_2(a,b)\phi^b(a)\psi(b) &=
  \phi(a){}^a\psi(b)\alpha_2(a,b) \\ G^{001} & =
  G^{a,b,\bar{b}\bar{a}m}_{ab,m,\bar{a}m}: & \tilde{\alpha}_3 \theta &=
  \alpha_3 \lc \phi\\ G^{011} &=
  G^{a,m,\bar{m}\bar{a}b}_{am,b,\bar{a}b}: & \tilde{\beta}_1(a,b)
  {}^{\bar{a}}\phi(a) {}^{\bar{a}}\omega(b) &= \omega(\bar{a}b)
  \theta(a,\bar{a}b) \beta_1(a,b) \\ G^{101} &=
  G^{m,a,\bar{a}\bar{m}b}_{ma,b,\bar{m}b}: & \tilde{\beta}_2(a,b)
  \psi^{\bar{a}}(a) \omega^{\bar{a}}(b) &= {}^b\bar{\phi}(a) \omega(b)
  \beta_2(a,b) \\ G^{110} &= G^{b\bar{a}\bar{m},m,a}_{b\bar{a},b,ma}: &
  \tilde{\beta}_3(a,b) \theta(b\bar{a},a)
  \bar{\omega}^{b\bar{a}}(b\bar{a}) &= \bar{\omega}^{b\bar{a}}(b)
  \psi^{\bar{a}}(a) \beta_3(a,b) \\ G^{111} &=
  G^{m\bar{a},a\bar{m}b,\bar{b}m}_{b,m,a}: & \tilde{\gamma}(a,b)
  \omega^a(b) \phi(b) &= {}^b\bar{\omega}^a(a) \psi(a) \gamma(a,b)
  \maybreak{}\\ P^{0000} &= P^{a,b,c,d}_{ab,abc,abcd,bcd,cd}: & 1 &=
  \delta \alpha \\ P^{0001} &=
  P^{a,b,c,\bar{c}\bar{b}\bar{a}m}_{ab,abc,m,\bar{a}m,\bar{b}\bar{a}m}:
  & 1 &= \alpha \lc \alpha_3 \\ P^{1000} &=
  P^{m\bar{c}\bar{b}\bar{a},a,b,c}_{m\bar{c}\bar{b},m\bar{c},m,abc,bc}:
  & \rc \alpha_1 &= \alpha \\ P^{0010} &=
  P^{a,b,\bar{b}\bar{a}m\bar{c},c}_{ab,m\bar{c},m,\bar{a}m,\bar{b}\bar{a}m}:
  & \alpha_3(a,b) \alpha_2(ab,c) &= {}^a\alpha_2(b,c) \alpha_2(a,c)
  \alpha_3^c(a,b) \\ P^{0100} &=
  P^{a,\bar{a}m\bar{c}\bar{b},b,c}_{m\bar{c}\bar{b},m\bar{c},m,\bar{a}m,bc}:
  & \alpha_1(b,c) \alpha_2(a,bc) &= \alpha_2^c(a,b) \alpha_2(a,c)
  {}^a\alpha_1(b,c) \maybreak{}\\ P^{0011} &=
  P^{a,b,m,\bar{m}\bar{b}\bar{a}c}_{ab,abm,c,\bar{a}c,\bar{b}\bar{a}c}:
  & \alpha(a,b,\bar{b}\bar{a}c) \beta_1(ab,c) &=
  {}^{\ov{ab}}\alpha_3(a,b) {}^{\bar{b}}\beta_1(a,c) \beta_1(b,\bar{a}c)
  \\ P^{1100} &=
  P^{c\bar{b}\bar{a}\bar{m},m,a,b}_{c\bar{b}\bar{a},c\bar{b},c,mab,ab}:
  & \alpha(c\bar{b}\bar{a},a,b) \beta_3(ab,c) &= \beta_3(a,c\bar{b})
  \beta_3^{\bar{a}}(b,c) \alpha_1^{\ov{ab}}(a,b) \\ P^{0101} &=
  P^{a,m,b,\bar{b}\bar{m}\bar{a}c}_{am,amb,c,\bar{a}c,\bar{m}\bar{a}c}:
  & \beta_1(a,c) {}^{\bar{a}}\beta_2(b,c) &= \beta_2(b,\bar{a}c)
  \beta_1^{\bar{b}}(a,c) {}^{\bar{a}}\alpha_2^{\bar{b}}(a,b) \\ P^{1010}
  &=
  P^{c\bar{b}\bar{m}\bar{a},a,m,b}_{c\bar{b}\bar{m},c\bar{b},c,amb,mb}:
  & \ov{{}^{b\bar{c}}\beta_2^a}(a,c) \beta_3(b,c) &=
  {}^{\bar{a}}\alpha_2^{\bar{b}}(a,b) {}^{\bar{a}}\beta_3(b,c)
  \ov{{}^{b\bar{c}}\beta_2^a}(a,c\bar{b}) \\ P^{0110} &=
  P^{a,m,\bar{m}\bar{a}c\bar{b},b}_{am,c\bar{b},c,\bar{a}c,\bar{m}\bar{a}c}:
  & \beta_1(a,c) \ov{\beta_3^{b\bar{c}a}}(b,c) &=
  \ov{\beta_3^{b\bar{c}a}}(b,\bar{a}c) \alpha(a,\bar{a}c\bar{b},b)
  \beta_1(a,c\bar{b}) \\ P^{1001} &=
  P^{m,a,b,\bar{b}\bar{a}\bar{m}c}_{ma,mab,c,\bar{m}c,\bar{a}\bar{m}c}:
  & \beta_2(a,c) \beta_2^{\bar{a}}(b,c) &= {}^c\bar{\alpha}_3(a,b)
  \beta_2(ab,c) \alpha_1^{\ov{ab}}(a,b) \maybreak{}\\ P^{0111} &=
  P^{a,\bar{a}m\bar{c},c\bar{m}b,\bar{b}m}_{m\bar{c},b,m,\bar{a}m,c}: &
  \alpha_2(a,c) \gamma(c,b) &= {}^a\gamma(c,\bar{a}b)
  \alpha_3(a,\bar{a}b) {}^a\beta_1^c(a,b) \\ P^{1110} &=
  P^{m\bar{c},c\bar{m}b,\bar{b}m\bar{a},a}_{b,m\bar{a},m,c,\bar{b}m}: &
  \alpha_2(b,a) \gamma(c,b) &= {}^b\beta_3^a(a,c) \alpha_1(c\bar{a},a)
  \gamma^a(c\bar{a},b) \\ P^{1011} &=
  P^{m\bar{c}\bar{a},a,c\bar{m}b,\bar{b}m}_{m\bar{c},b,m,ac,c}: &
  \alpha_1(a,c) \gamma(c,b) &= \beta_2^{ac}(a,b) \gamma(ac,b)
  {}^b\bar{\beta}_1 {}^c(a,ac) \\ P^{1101} &=
  P^{m\bar{c},c\bar{m}b,a,\bar{a}\bar{b}m}_{b,ba,m,c,\bar{b}m}: &
  \gamma(c,b) \alpha_3(b,a) &= {}^b\bar{\beta}_2{}^c(a,c) \gamma(c,ba)
  {}^b\bar{\beta}_3{}^c(a,ba) \\
  P^{1111} &=
  P^{a\bar{m},m,\bar{m}e,\bar{e}m\bar{a}b}_{a,a\bar{m}e,b,m\bar{a}b,d}:
  & \delta_{b,ad} {}^e\bar{\beta}_1(a,b) \beta_3(d,b) &= \sum_{c \in eA}
  {}^e\bar{\gamma}^a(c,a) \bar{\beta}_2^a(c,b) \gamma^{\bar{d}}(d,c)
\end{align*}
where $P^{1111}$ assumes $\bar{b}ad \in A$.

\begin{defn}
  $f$ is \emph{normal} if $\beta_1(-,1) = \beta_2(-,1) \equiv 1$.
\end{defn}

\begin{lem}
  Every fusion system on $L$ is gauge equivalent to a normal one.
\end{lem}

\begin{proof}
  Given arbitrary $f$, we construct $\xi$ making $\tilde{f}$ normal.  By
  $G^{011}$,
  \[ \tilde{\beta}_1(a,1) = \frac{ \omega(\bar{a}) \theta(a,\bar{a})
    \beta_1(a,1) }{ {}^{\bar{a}}\omega(1) {}^{\bar{a}}\phi(a) }\] so we
  can choose $\omega$ making $\tilde{\beta}_1(a,1) = 1$.  By $G^{101}$,
  \[ \tilde{\beta}_2(a,1) = \frac{ \bar{\phi}(a) \omega(1) \beta_2(a,1)}
  {\psi^{\bar{a}}(a) \omega^{\bar{a}}(1)}\] so we can choose $\psi$
  making $\tilde{\beta}_2 = 1$.  Then any choice of normalized $\phi,
  \theta$ makes $\tilde{f}$ a fusion system on $L$ by
  Lemma~\ref{lem-gaugify}.
\end{proof}

\begin{lem}
  $\Psi$ restricts to an isomorphism on the category of normal fusion
  systems on $L$.
\end{lem}

\begin{proof}
  It is routine to check $\mathcal{Y}$ and $\mathcal{X}$ are categories.
  Let $\mathcal{Z}$ be the category of normal fusion systems, a full
  subcategory of $\mathcal{Y}$, and let $\Xi = \Psi|_{\mathcal{Z}}$.
  Let $f \in \obj \mathcal{Z}$, and $(\chi,\ups,\tau) = \Xi f =
  (\alpha_2,\alpha_3,\gamma(1,1))$. Then
  \begin{align*}
    \alpha &= \frac{1}{\lc \ups} \colon M^3 \to \F^\times &(P^{0001})\\
    \beta_1(b,a) &= \frac{ \alpha(\bar{a},b,\bar{b}a) }{
      {}^{\bar{b}a}\ups(\bar{a},b)}
    &(P^{0011})\\
    \beta_2(b,a) &= \frac{\ups(a,\bar{a})}{\ups^{\bar{b}}(a,\bar{a})} \chi^{\bar{b}}(a,b) &(P^{0101})\\
    \alpha_1(b,c) &= \frac{1}{\bar{\ups}^{bc}(b,c)}     &(P^{1001},P^{0100})\\
    \chi &\text{ is an $\ups$-biderivation} &(P^{0010}, P^{0100})
    \maybreak{}\\
    \gamma(a,1) &= \tau \bar{\ups}(\bar{a},a) &(P^{1011})\\
    \gamma(c,a) &= \lexp{a}{\left( \frac{ \tau \bar{\ups}(\bar{c},c) }{ \ups^c(\bar{a},a) } \right)} \frac{1}{\chi(a,c)} &(P^{0111})\\
    \tau &\in B^\times     &(\text{invertibility}) \maybreak{}\\
    \beta_3(a,c) &= \frac{ \bar{\ups}(a,\bar{c}) \tau^{\bar{a}} }{ \alpha(a,\bar{c},c) \alpha(a\bar{c},c\bar{a},a)\tau} &(P^{1110})\\
    \tau \bar{\tau} &\in \F^\times &(P^{1101})\\
    \chi &\text{ is $\tau$-quasisymmetric} &(P^{1010})\\
    \tau \bar{\tau} &\equiv \nicefrac{1}{\abs{A}} &(P^{1111})\\
    \chi|_{A \times A} &\text{ is nondegenerate} &(P^{1111})
  \end{align*}
  Thus $\Xi$ is well-defined and injective on the level of objects.

  To check $\Xi$ is surjective on the level of objects, suppose
  $(\chi,\ups,\tau) \in \mathcal{X}$.  We must show $f$ defined via the
  above equations, with $\alpha_2 = \chi$ and $\alpha_3 = \ups$, is in
  $\obj \mathcal{Z}$.  It suffices to show $f$ is a fusion system on
  $L$.  By construction $f$ satisfies admissibility.  To check
  invertibility, let $p,q \in S$.  Since $\tau \bar{\tau} \equiv
  \nicefrac{1}{\abs{A}}$ and $\chi|_{A \times A}$ is a nondegenerate
  symmetric bicharacter, the matrix
  $(\nicefrac{{}^q\tau}{\chi(y,x)})_{x,y \in A}$ has inverse
  $({}^q\bar{\tau}\chi(y,x))_{x,y \in A}$.  Since
  \[ \gamma(px,yq) = \frac{ {}^q\bar{\ups}(\bar{x}\bar{p},px) }{
    \chi(q,px) } \frac{{}^q\tau}{\chi(y,x)} \frac{ \ups(y,q) }{
    \chi(y,p) \ups^p(y,q) {}^q\ups^p(\bar{q}\bar{y},yq) }\] for $x,y \in
  A$, the matrix $(\gamma(px,yq))_{x,y \in A}$ is invertible over $B$.
  Equivalently, for all $m \in M$ the recoupling matrix
  $F^{m\bar{p},p\bar{m}q,\bar{q}m}_m = (\gamma(u,v)(m))_{u \in pA, v \in
    qA}$ is invertible over $\F$.  Therefore $f$ satisfies the
  invertibility axiom by Lemma~\ref{lem:cosetmats}.  It is routine to
  verify the pentagon axiom, and trivial to verify the triangle axiom.
  Rigidity means $\alpha(a,\bar{a},a)\alpha(\bar{a},a,\bar{a}) = 1$ and
  $\tau(m) = ((F^{\bar{m},m,\bar{m}}_{\bar{m}})^{-1})_{1,1}$ for all $a
  \in S$ and $m \in M$.  The first condition follows from $\delta
  \alpha(a,\bar{a},a,\bar{a}) = 1$.  The second follows from
  \begin{align*}
    F^{\bar{m},m,\bar{m}}_{\bar{m}} &= \left( \frac{\bar{\ups}(\bar{x},x) \tau}{ \chi(y,x) \ups(\bar{y},y) } (\bar{m}) \right)_{x,y \in A}\\
    (F^{\bar{m},m,\bar{m}}_{\bar{m}})^{-1} &= \left(
      \frac{\ups(\bar{x},x) \bar{\tau}\chi(y,x) }{ \bar{\ups}(\bar{y},y)
      } (\bar{m})\right)_{x,y \in A}
  \end{align*}
  as $\chi$ is normalized and $\tau(m) = \bar{\tau}(\bar{m})$.  Thus $f$
  is a fusion system on $L$, showing $\Xi$ is surjective on the level of
  objects.

  Suppose $\xi \in \mor_{\mathcal{Z}}(f,\tilde{f})$.  Letting
  $(\tilde{\chi},\tilde{\ups},\tilde{\tau}) = \Xi \tilde{f}$, we have
  \begin{align*}
    \om(a) &= \frac{ \lexp{a}{\phi}(\bar{a}) {}^a\varsigma }{ \theta(\bar{a},a) } &(G^{011})\\
    \psi(a) &= \frac{\bar{\phi}^a(a) \varsigma^a }{ \varsigma} &
    (G^{110}) \maybreak{}\\
    \frac{\tilde{\ups} }{\ups} &= \frac{\lc \phi}{\theta}  &(G^{001})\\
    \frac{\tilde{\chi}(a,b)}{\chi(a,b)} &=
    \frac{\phi(a){}^a\psi(b)}{\phi^b(a)\psi(b)} = \frac{ \phi(a)
      {}^a\bar{\phi}^b(b) \lexp{a}{\varsigma}^b \cdot \varsigma }
    { \phi^b(a) \bar{\phi}^b(b) \lexp{a}{\varsigma} \varsigma^b}  &(G^{010})\\
    \frac{\tilde{\tau}}{\tau} &= \frac{\bar{\varsigma}}{\varsigma} &(G^{111})\\
  \end{align*}
  Since $\theta \colon M \to \F^\times$, we see $\Xi$ is well-defined
  and injective on each morphism space.  It is routine to check $\Xi$ is
  surjective on each morphism space.  Therefore $\Xi$ is an isomorphism.
\end{proof}

\begin{proof}[Proof of Theorem~\ref{thmKludge-main}]
  The two preceding lemmas suffice.
\end{proof}

\section*{Acknowledgments}
It is a pleasure to thank my thesis advisor, Zhenghan Wang, for
essential background and advice; Parsa Bonderson and Joost Slingerland
for vital, generous feedback; and Tobias Hagge and Scott Morrison for
useful discussions.

\bibliography{GenTYbib} \bibliographystyle{hplain} \nocite{*}
\end{document}